\documentclass[12pt]{article}
\usepackage[utf8]{inputenc}
\usepackage{graphicx}
\graphicspath{ {images/} }
\usepackage{amsmath, amssymb}
\usepackage{geometry}
\usepackage{float}  
\usepackage{graphicx}
\usepackage{cite}
\usepackage{hyperref}
\usepackage{lmodern}
\usepackage{placeins} 
\usepackage{tikz}
\usepackage{pgfplots}
\pgfplotsset{compat=1.18}

\begin{document}

\begin{titlepage}
\vspace*{\fill}
\begin{center}
{\LARGE \textbf{Real-Time Inviscid Fluid Dynamics and Aero-acoustics on a Sphere} \par}
\vspace{1.5cm}

{\large Madhusraba Sinha and Jan Stratmann \par}
\vspace{0.5cm}

{\normalsize Department of Mathematics, TU Dortmund University \par}
{\normalsize \texttt{madhusraba.sinha@tu-dortmund.de} \par}
{\normalsize \texttt{jan3.stratmann@tu-dortmund.de} \par}
\end{center}
\vspace*{\fill}
\end{titlepage}

\newpage

\begin{titlepage}
\vspace*{\fill}
\begin{center}
{\Large \bfseries Abstract \par}
\vspace{0.5cm}

\begin{minipage}{0.9\textwidth}
\small
Real-time fluid and aeroacoustic simulation on complex surfaces can have interactive applications - from globe-based weather visualizations to immersive computer games with physically accurate wind and sound. However, conventional grid-based solvers struggle with numerical instability near surface singularities, and mesh-based approaches lack a straightforward path to solving partial differential equations (PDEs) with stable, high-order accuracy.

Our model presents a unified framework for real-time inviscid fluid simulation and aeroacoustics on spherical surfaces with embedded obstacles, combining the Closest Point Method (CPM), projection-based Navier–Stokes solvers, and the Ffowcs Williams–Hawkings (FWH) analogy. CPM enables surface PDEs to be solved in a Cartesian embedding without parametrization by restricting computation to a narrow band around the sphere. Each band point is mapped to its nearest surface location, where band operators project results onto the local tangent space. Surface obstacles are modelled with signed distance functions (SDFs), enforcing no-slip velocity constraints and Bernoulli-based pressure adjustments for consistent real-world boundary interactions. Aeroacoustic sources are computed directly from surface pressure force derivatives and mapped to real-time audio via frequency and amplitude modulation with artifact-suppressing hysteresis smoothing.

Our findings from this model simulate the behaviour of inviscid fluid on spherical surfaces while generating sound using the pressure of the fluid flowing on the surface. This approach gives results that offer stability, geometric consistency, and support applications in scientific visualization, virtual reality, and educational tools.

All source code is made publicly available via the GitHub repository https://github.com/jstr01/RealTime-CFD-Acoustics-Sphere.
\newline


\end{minipage}
\end{center}
\vspace*{\fill}
\end{titlepage}

\section{Introduction}

Simulating fluid flow and sound on curved surfaces can be used in many areas, from weather prediction to virtual environments. These simulations make it possible to create realistic audio-visual effects in various simulations and graphics in computer games. But once we move away from flat Cartesian geometry, the usual numerical tools have trouble handling the more complex shapes. Spherical domains introduce coordinate singularities, uneven grid spacing, and metric complexity — each of which can destabilize simulations or drive computational costs beyond practical limits. This work presents a novel framework that combines advanced numerical techniques to overcome these challenges, which enables stable, real-time fluid and aeroacoustic simulations directly on spherical surfaces with obstacles placed on it. Our approach merges real-time stable fluid dynamics with surface-bound partial differential equation (PDE) solvers and efficient acoustic modelling. This fusion retains physical plausibility while remaining computationally lightweight. 

Stam introduced the Stable Fluids method in \cite{Stam1999}, shaping real-time fluid simulation by utilising semi-Lagrangian advection and implicit diffusion to achieve unconditional numerical stability at large time steps. While Stam’s model serves as an excellent foundation for interactive applications, it is mainly designed for flat Euclidean domains and does not apply to curved surfaces with introduced obstacles. Building on the fractional-step incompressibility enforcement introduced by Chorin in \cite{Chorin1968} and developed by Brown, Cortez, and Minion in \cite{Brown2001}, our work extends these projection methods to handle Neumann boundary conditions defined implicitly by signed distance functions (SDFs). We combine Stam’s stable semi-Lagrangian method with projection schemes that work with surface geometry, which lets us run stable fluid simulations on curved surfaces while accounting for obstacles, without the need for remeshing or parameterization. Previously, Stam introduced Catmull–Clark subdivision surfaces in \cite{StamSubdivision1998}, which create smooth geometry through recursive refinement. These surfaces are common in graphics, but they are designed for geometry and don’t directly support PDE solvers. Our framework integrates subdivision surfaces with the Closest Point Method (CPM), allowing numerical operators to act consistently on the smooth geometry and the physical simulation domain. Osher and Sethian introduced level-set and SDF methods in \cite{Osher1988}, which use flexible implicit obstacle representations. These methods work smoothly with fluid solvers and are clearly defined in numerical terms. Yang, Corse, Lu, Wolper, and Jiang further developed spherical fluid simulation using latitude-longitude grids in \cite{YangMay2019} but encountered pole singularities and grid non-uniformity, which reduced numerical accuracy. The Closest Point Method (CPM) developed by Ruuth and Merriman in \cite{Ruuth2008} and generalized by Macdonald in \cite{MacdonaldThesis2008} solves PDEs on surfaces by using them in Euclidean space and applying Cartesian finite differences within a narrow band. Our solver addresses the singularity and interpolation issues of initial spherical solvers by utilising CPM’s strategy combined with sparse grid acceleration, ensuring efficient, stable, and accurate fluid simulation on spherical domains. Sparse grid techniques, reviewed by Bungartz and Griebel in \cite{Bungartz2004}, reduce the problem of dimensionality in high-resolution simulations by applying hierarchical basis functions. We use sparse grids to handle interpolation between the simulation domain and the narrow band around the surface. This setup keeps accuracy and allows real-time performance in fluid-acoustic simulations. The foundational aeroacoustic model by Ffowcs Williams and Hawkings in \cite{FWH1969} represents fluid-generated sound as a wave equation with turbulence and surface motion as sources. Further reviews by Brentner and Farassat in \cite{Brentner2003} and studies by Shen and Sørensen in \cite{Shen26May1999} validate the feasibility of aeroacoustic modelling coupled with incompressible flow solvers. We embed aeroacoustic source computation directly on CPM surfaces while keeping the obstacle geometry consistent. This approach generates sound predictions that align with the fluid simulation, and their accuracy is supported by the convergence results shown in Section 5.

Our unified framework leverages subdivision surfaces, SDF obstacles, CPM-based PDE solvers, and sparse grid interpolation to keep surface representations consistent across fluid, boundary conditions, and aeroacoustic calculations. This tight integration improves both numerical stability and physical realism. This paper details how we developed this framework, showing its solid mathematical basis, efficient computation, and ability to simulate flows on various geometries, interactions with complex obstacles, and real-time coupled sound generation.

The next section discusses the foundations of CPM embedding and the reasoning behind this choice. Section 3 provides the governing equations and computational fluid dynamics methods used to develop the framework for this model. To understand the aeroacoustics and their relevance to this work, Section 4 presents the aeroacoustic model and its numerical implementation using the FW-H analogy. Section 5 presents the results, which explains the accuracy and efficiency of our model. Finally, Section 6 concludes with a discussion and future research.

\section{Methodological Comparison for Surface PDEs}

In this section, we outline the specifications we want our method to achieve for real-time fluid and acoustic simulations on a sphere. We would need a framework that can handle complex spherical geometry without running into problems like coordinate singularities, is easy to implement using standard numerical operations on a uniform grid, and can support real-time performance with efficient data handling. At the same time, it should connect naturally with physical models such as fluid flow, pressure, velocity, and acoustic fields, and remain flexible enough to handle more complex scenarios like turbulent flows or coupled sound propagation. The subsections that follow go through different approaches and the comparative analysis of those approaches. 

\subsection{Catmull–Clark Subdivision Surfaces}
Originally designed for graphical smoothness, Catmull-Clark subdivision surfaces recursively subdivide coarse polygonal meshes into smooth, continuous surfaces. It would require watertight meshes. Adding or deforming obstacles forces global re-meshing, which breaks real-time performance. While they are suitable for representing geometry, these processes do not inherently support the breakdown of scalar or vector fields or the efficient computation of PDE operators needed for surface-bound fluids and acoustics. To handle derivatives, interpolate velocity fields, or apply conservation laws on a surface made from these methods, we need additional frameworks. This often involves projecting or parameterizing the mesh by introducing metric tensors into every PDE, which increases computational complexity and memory usage than that required for CPM.  

\subsection{Staggered Grids on the Sphere}
Staggered grid methods store scalar and vector values at offset spatial locations to reduce pressure and velocity decoupling. This grid staggering is common in finite volume or finite difference CFD solvers. Staggered grid setups in a sphere usually follow latitude-longitude (spherical) or cubed-sphere patterns. Staggered grids assume uniform Cartesian volumes, but on spheres it results in breaking mass conservation and requires correction terms. In this case, obstacles would rarely align with the grid which in turn would produce “staircase artifacts” in pressure and velocity fields. Spherical coordinates work well in flat shapes, but they have issues at the poles and with uneven grid cell spacing. This can lead to stability problems and uneven resolution, especially in flows dominated by advection. Also, it is inconvenient to apply physically consistent boundary conditions and divergence-free conditions with a non-uniform spherical grid.

\subsection{Closest Point Method (CPM)}
The Closest Point Method (CPM) provides a simple and flexible way to solve PDEs on surfaces. It works by extending the solution off the surface through closest-point projection, which allows the use of standard Cartesian finite difference operators in 3D without needing complex meshes, parameterizations, or coordinate mappings. In our implementation, we combine this approach with a narrow-band masking technique around the sphere, so computations are focused only where they are needed, accurately representing the surface geometry without handling the full 3D domain. To solve the resulting linear systems efficiently, we use a Cholesky decomposition, which keeps the computations stable and fast. CPM is easy to implement, works with time-stepping and conservative schemes, and can handle a wide range of PDEs, including diffusion, advection, and wave equations. It also naturally manages arbitrary surfaces and obstacles through implicit representations like signed distance functions (SDFs), works well with more complex models because it handles coupled PDEs consistently, and can be extended to other manifolds beyond spheres. Overall, CPM provides a computationally efficient, stable, and versatile framework that aligns well with our goals for real-time fluid and acoustic simulations.

\section{Governing Equations and Numerical Methods}

This section describes the mathematical and computational framework used to simulate incompressible fluid flow around a sphere with surface obstacles. It also covers the generation of aeroacoustic sound as a result. The main setup connects a fluid solver using the Closest Point Method (CPM) with an aeroacoustic solver that follows the Ffowcs Williams-Hawkings (FW-H) analogy. The results are turned into sound in real time.

\subsection{Fluid Dynamics Model}
The motion of an incompressible, inviscid fluid is described by the Euler equations, which represent the conservation of momentum and mass. The momentum equation captures the main inertial and pressure forces responsible for sound generation at low Mach numbers,
\begin{equation}
\frac{\partial \mathbf{u}}{\partial t} + (\mathbf{u} \cdot \nabla)\mathbf{u}
= -\frac{1}{\rho_0} \nabla p + \mathbf{f}
\label{eq:euler}
\end{equation}
while the continuity equation enforces incompressibility:
\begin{equation}
\nabla \cdot \mathbf{u} = 0,
\label{eq:incompressibility}
\end{equation}
In our simulations, we use the inviscid approximation ($\nu = 0$), focusing on dynamics driven primarily by inertia and pressure. A localized forcing term $\mathbf{f}$ is applied near the top of the sphere to initiate and sustain motion. Boundary conditions are imposed on obstacles as described in the Subsection Modelling Surface Obstacles, and average value conditions are applied wherever necessary to maintain global flow consistency.

\subsection{Closest Point Method for Surface PDEs}
Simulating PDEs directly on curved surfaces is difficult. We use the Closest Point Method (CPM) \cite{MacdonaldThesis2008,Ruuth2008}. This method embeds the PDE into a narrow volumetric band around the surface. It allows us to apply standard Cartesian finite-difference methods directly.

For any point $\mathbf{x}$ in the band $\Omega$, the surface field is extended by assigning the value at its closest point on the sphere:
\begin{equation}
u(\mathbf{x}) = u(\text{CP}(\mathbf{x})), \quad
\text{CP}(\mathbf{x}) = R \frac{\mathbf{x}}{\|\mathbf{x}\|},
\end{equation}
where $R$ is the sphere radius.  

This narrow-band embedding allows for advection, pressure projection, and interpolation back to the surface using the closest-point mapping.

\subsection{Numerical Discretization and Projection Method}
Equations \eqref{eq:euler}–\eqref{eq:incompressibility} are solved using a projection method as used by Chorin in \cite{Chorin1968}, which separates the velocity advection from the pressure projection step. First, an intermediate velocity $\mathbf{u}^*$ is calculated while ignoring pressure effects. This involves applying advection and any external forcing, using a semi-Lagrangian backtracking scheme to maintain numerical stability:
\begin{equation}
\frac{\mathbf{u}^* - \mathbf{u}^n}{\Delta t}
= -(\mathbf{u}^n \cdot \nabla)\mathbf{u}^n + \mathbf{f}.
\end{equation}

Next, the intermediate velocity is corrected by solving a Poisson equation for pressure, ensuring incompressibility, and then updating the velocity field accordingly:
\begin{equation}
\nabla^2 p^{n+1} = \frac{\rho_0}{\Delta t} \nabla \cdot \mathbf{u}^*,
\label{eq:poisson}
\end{equation}
\begin{equation}
\mathbf{u}^{n+1} = \mathbf{u}^* - \frac{\Delta t}{\rho_0}\nabla p^{n+1}.
\end{equation}

The pressure solve represents the main computational bottleneck of the method. In our implementation, we construct sparse Laplacian operators with Neumann boundary conditions applied on the narrow band around the surface. For small grids, Cholesky factorization is precomputed to speed up the solve, while larger systems are handled using a Conjugate Gradient iterative method. This approach balances efficiency and accuracy, enabling real-time or near real-time updates of the velocity field on the sphere.

\subsection{Modelling Surface Obstacles}
Surface obstacles are modeled as rigid, immovable bodies using Signed Distance Functions (SDFs). Each obstacle is defined by its centre, radius, and height. The obstacles affect the simulation in two main ways. Inside an obstacle, the velocity is set to zero to represent a no-slip condition, while on the surface of the obstacle, the velocity is projected onto the tangent plane to enforce no-penetration. For visualization, the sphere mesh is modified with Gaussian-shaped bumps at the locations of the obstacles. This approach ensures that the fluid interacts physically consistently with the embedded surface features while maintaining an accurate geometric representation.

\subsection{Software Implementation}
The solver is built in Python. It uses \texttt{NumPy} and \texttt{SciPy} for the core numerical computations. Sparse linear systems, such as the ones from the pressure Poisson equation, are handled with precomputed Cholesky factorization for smaller grids and Conjugate Gradient iteration for larger ones. Advection is handled using a semi-Lagrangian scheme, while incompressibility is enforced through a pressure projection step. The Closest Point Method (CPM) and narrow-band masking allows us to efficiently solve surface PDEs directly on the sphere, while obstacles are represented using signed distance functions (SDFs). For visualization, we use \texttt{PyVista} within a \texttt{PyQt6} GUI, which makes it easy to interactively explore the sphere and the obstacles, and real-time audio is generated with \texttt{sounddevice}. We explain more about the real-time audio generation in the next section. The whole system is designed modularly, separating fluid simulation, acoustics, and rendering. This keeps the solver efficient while also making it flexible and easy to maintain, such that alternative solvers, more complex physics, or updated visualization and audio modules can be added without much difficulty.

\section{Aeroacoustic Modelling and Sound Synthesis}

\subsection{Layout}
The prediction of sound from unsteady flow is a significant challenge in computational aeroacoustics (CAA) because it involves modelling the complex fluid motion while also accounting for how these changes produce acoustic waves. To address this issue, we use a hybrid framework that merges fluid simulation and real-time sound synthesis. We simulate the fluid dynamics on the surface of a sphere using the Closest Point Method (CPM), and we obtain the radiated sound through the Ffowcs Williams-Hawkings (FW-H) acoustic analogy. We then synthesize the resulting pressure fluctuations into real-time audio. This allows for both physically-based far-field predictions and an auditory representation of flow changes that can be observed interactively. By combining CPM-based surface flow modelling with FW-H acoustics, our method effectively predicts far-field pressure variations. Additionally, frequency-domain analysis and additive synthesis enhance the sound's richness, turning complex fluid-acoustic interactions into an engaging audio experience.

\subsection{Theoretical Basis: FW–H Acoustic Analogy}
The acoustic radiation is modelled using the Ffowcs Williams-Hawkings (FW-H) equation in \cite{FWH1969}. This equation reformulates the compressible Navier–Stokes equations into an inhomogeneous wave equation, which makes it easier to describe how flow fluctuations turn into sound. In the fluid simulation we treat the flow as incompressible, which simplifies the problem by focusing on the main dynamics without having to resolve for small density changes. However, for sound, compressibility has to be included because acoustic waves are mainly pressure fluctuations that travel through the medium. The FW-H equation provides a way to connect these two viewpoints, letting us use an incompressible flow solver while still modelling the physics of sound radiation. For a rigid and stationary surface $S$, the far-field pressure fluctuation $p'(\mathbf{x},t)$ at the observer location $\mathbf{x}$ is given by:

\begin{equation}
p'(\mathbf{x},t) \;=\; \frac{1}{4\pi c_0} 
\frac{\partial}{\partial t} \int_{S} 
\left[\frac{p(\mathbf{y},\tau^*) \, \mathbf{n}(\mathbf{y}) \cdot \mathbf{i}_r}{r}\right] dS ,
\end{equation}

where:
\begin{itemize}
    \item $c_0$ is the ambient speed of sound,
    \item $p(\mathbf{y},t)$ is the unsteady pressure at surface point $\mathbf{y} \in S$,
    \item $\mathbf{n}(\mathbf{y})$ is the unit outward normal at $\mathbf{y}$,
    \item $r = \|\mathbf{x}-\mathbf{y}\|$ is the observer distance,
    \item $\mathbf{i}_r = (\mathbf{x}-\mathbf{y})/r$ is the radiation direction,
    \item $\tau^* = t-r/c_0$ is the retarded emission time.
\end{itemize}

For compact, low-Mach sources, the FW–H expression can be simplified. A source is considered compact when its size is much smaller than the wavelength of the sound it generates, which means it can be treated like a single point. Low-Mach refers to flows with a Mach number below one, where the flow speed is below the speed of sound. In this case, non-linear compressibility effects are minimal, and the main contribution to the radiated sound comes from unsteady pressure forces acting on the surface. Under these assumptions, the FW–H equation reduces to:

\begin{equation}
p'(\mathbf{x},t) \;\approx\; \frac{1}{4\pi c_0 r} \, 
\mathbf{i}_r \cdot \frac{d\mathbf{F}}{dt}(t),
\end{equation}

where the unsteady aerodynamic force acting on the body is

\begin{equation}
\mathbf{F}(t) \;=\; \int_{S} p(\mathbf{y},t) \, \mathbf{n}(\mathbf{y}) \, dS .
\end{equation}

So, the far-field sound is directly proportional to the time derivative of the aerodynamic force acting on the body.

\subsection{Frequency Analysis for Timbre Enrichment}
While the FW-H formulation produces a monopole-like signal, we can gain more detailed timbral information by looking at the frequency content of $\tfrac{d\mathbf{F}}{dt}(t)$. We apply a sliding temporal window to this signal and perform a Fourier transform:

\begin{equation}
\mathbf{\hat{F}}(\omega) \;=\; \int_{-\infty}^{\infty} 
\frac{d\mathbf{F}}{dt}(t) \, e^{-i \omega t} \, dt ,
\end{equation}

which gives us the spectral amplitudes and phases of the fluctuating force. Peaks in $|\mathbf{\hat{F}}(\omega)|$ relate to coherent vortical structures and periodic shedding phenomena in the flow.

We identify the dominant frequencies $\{\omega_k\}$ by finding the local maxima of $|\mathbf{\hat{F}}(\omega)|$. The associated magnitudes $|\mathbf{\hat{F}}(\omega_k)|$ serve as weights for the next stage of synthesis.

While the current monopole-based FW-H model reproduces the global force fluctuations, extending the model to include dipole and quadrupole enables a more physically complete description of the aerodynamic sound generation. The model would then transform to {Lighthill's acoustic analogy:

\begin{equation}
\nabla^2 p' - \frac{1}{c_0^2}\frac{\partial^2 p'}{\partial t^2} 
= -\nabla \cdot \left( \nabla \cdot \mathbf{T} \right),
\end{equation}
where the Lighthill stress tensor is defined as
\begin{equation}
T_{ij} = \rho u_i u_j + (p' - c_0^2 \rho') \delta_{ij} - \tau_{ij}.
\end{equation}

The three primary acoustic sources that emerge under the assumption of incompressible flow are:
\begin{align*}
\text{Monopole (thickness noise):} & \quad \frac{\partial^2 \rho}{\partial t^2},\\
\text{Dipole (loading noise):} & \quad \nabla \cdot \mathbf{F},\\
\text{Quadrupole (shear noise):} & \quad \nabla \cdot \big( \nabla \cdot (\rho \mathbf{u} \otimes \mathbf{u}) \big).
\end{align*}

In the case of the extended model, dipoles represent unsteady forces on solid surfaces, and quadrupoles describe turbulence and vortex interactions in the fluid volume. Including these terms would increase both spatial and temporal resolutions. The model would evolve from a surface-only integration (\( \mathcal{O}(N_\text{surface}) \)) to a combined surface and volume computation (\( \mathcal{O}(N_\text{surface} + N_\text{volume}) \)). 

\subsection{Real-Time Sound Synthesis}  
The sound synthesis engine converts aeroacoustic quantities into an audible waveform. The instantaneous acoustic pressure magnitude, $p'(t)$, is used to control the overall amplitude envelope of the sound. A non-linear mapping,  
\begin{equation}  
A(t) \;=\; \tanh\!\big(\alpha \, |p'(t)|\big) ,  
\end{equation}  
is applied to prevent clipping and to recreate the way humans perceive loudness saturation, where $\alpha$ is a scaling factor. To capture the tonal characteristics of the sound, dominant frequencies $\omega_k$ are extracted from a spectral analysis of the aerodynamic forces, which are then used to drive individual sinusoidal oscillators. The resulting synthesized signal is given by   
\begin{equation}  
s(t) \;=\; A(t) \sum_{k=1}^{N}  
\frac{|\mathbf{\hat{F}}(\omega_k)|}{\sum_j |\mathbf{\hat{F}}(\omega_j)|}  
\; \sin(\omega_k t + \phi_k) ,  
\end{equation}  
where $\phi_k$ are phase terms and the magnitudes are normalized to maintain energy balance. This approach produces a composite waveform whose tonal characteristics reflect the dominant vortical modes of the flow, while the overall loudness corresponds to the calculated acoustic pressure. In cases where there are no distinct frequency peaks, such as in laminar flow regimes, the model defaults to a single dominant frequency $\omega_d$, ensuring continuous audio output:  
\begin{equation}  
s(t) \;=\; A(t) \sin(\omega_d t + \phi_d).  
\end{equation}  

\subsection{Coupling with the Flow Solver and Sound Synthesis}  
In our model, the aeroacoustic pipeline connects to the surface flow solver, allowing fluid dynamics to directly influence the generated sound. The pressure field computed in the flow solver is mapped onto the sphere surface using the Closest Point Method (CPM), and the unsteady aerodynamic force on the surface and obstacles is calculated through surface quadrature:
\[
\mathbf{F}(t) = \int_S p(\mathbf{y},t)\,\mathbf{n}(\mathbf{y})\,dS(\mathbf{y}).
\]
The time derivative of this force, $\frac{d\mathbf{F}}{dt}$, drives the compact-source FW–H model to determine far-field pressure fluctuations, $p'(\mathbf{x},t)$, at observer locations. Dominant frequencies are identified in real time from the aerodynamic force using spectral analysis and are assigned to sinusoidal oscillators in the sound synthesis engine. The amplitude of these oscillators is modulated by the instantaneous pressure magnitude, producing a composite waveform that reflects both the tonal character and loudness of the flow. This integrated approach ensures that variations in the fluid structures produce a corresponding and physically consistent acoustic output.

\section{Manufactured solutions}

To confirm that the solver behaves correctly, we use the Method of Manufactured Solutions (MMS). We implement this by choosing analytic expressions for velocity, density, and pressure and then constructing forcing terms that make them satisfy the equations exactly. This gives us a reference solution, which makes it easy to see how much error the solver produces. All analytic expressions are explained in the python file \texttt{manufactured\_solution.py}, and the full verification procedure is implemented in the python file \texttt{test\_manufactured\_solution.py}.

\subsection{Analytic fields}
The manufactured fields are built using smooth trigonometric functions on the unit cube. These functions are easy to evaluate but have some variation to properly test the solver. The velocity, density, and pressure fields are defined as:

\begin{align}
u_1(x,y,z) &= \sin(\pi x)\cos(\pi y)\sin(\pi z), \\
u_2(x,y,z) &= -\cos(\pi x)\sin(\pi y)\sin(\pi z), \\
u_3(x,y,z) &= \sin(\pi x)\sin(\pi y)\cos(\pi z), \\
\rho(x,y,z) &= 1.0 + 0.1 \sin(\pi x)\sin(\pi y)\sin(\pi z), \\
p(x,y,z) &= 1.0 + 0.05 \cos(\pi x)\cos(\pi y)\cos(\pi z).
\end{align}

These fields are divergence-free and include three-dimensional variation to test how the solver handles advection, diffusion, and the interaction between different flow variables properly. Since the solver evaluates the same analytic expressions internally, we can compare its output to the manufactured fields directly, which makes it easy to check whether the numerical results match as they should.

\subsection{Forcing terms}
To ensure that the manufactured fields satisfy the governing equations, we compute the forcing terms analytically. Since viscosity is not considered, the forcing only balances the non-linear advection and pressure-gradient terms in the Euler equations. Specifically:

\begin{itemize}
\item \texttt{forcing\_f(...)} provides the momentum forcing term $\mathbf{f}$:
\[\mathbf{f} = - \Big[ (\mathbf{u} \cdot \nabla) \mathbf{u} + \frac{1}{\rho_0} \nabla p \Big],\]
so that the manufactured velocity field satisfies the modified momentum equation exactly.
\item \texttt{scalar\_source\_S(...)} provides the source term for density:  
\[S = - \nabla \cdot (\rho \mathbf{u}),\]  
ensuring that the analytic density field is an exact solution of the transport equation.
\end{itemize}

This setup allows the solver to reproduce the analytic fields precisely, meaning any numerical error arises only from discretization.

\subsection{Numerical experiment setup}
The verification tests follow a simple pattern. For each grid resolution, the code generates a uniform grid and initializes the solver with \texttt{use\_manufactured\_solution=True}. Obstacles can also be included, and MMS remains valid because the forcing terms are defined everywhere.

After one timestep, we compare the numerical fields against the analytic ones. The tests check for the accuracy of the velocity and density fields. This is done by calculating their L2 errors, which should be below 0.25 for velocity and 0.4 for density at grid resolutions of $11^3$, $13^3$, and $15^3$. This shows that the solver is correctly capturing the prescribed fields. We also evaluate the divergence of the velocity field to ensure it remains reasonable, even though incompressibility is not explicitly enforced; the maximum divergence should be less than 5.0 on all tested grids. Finally, we test how the solver behaves when the grid is refined by running it on finer grids with resolutions $11$, $13$, $15$, and $17$ and collecting the L2 errors for velocity, density, and divergence. These checks help confirm that the solver handles advection and forcing terms properly and that the numerical solution behaves consistently as the grid is refined.

\subsection{Convergence rates}
To measure how the solver improves with refining the resolution, we calculate experimental convergence rates using the L2 errors from resolutions $11$, $13$, $15$, and $17$. For each grid size the spacing is defined as $h = 1/(r-1)$, and we fit a line to $\log(\text{error})$ versus $\log(h)$. The slope of this line represents the rate at which the error decreases as the grid becomes finer.

Because the manufactured fields are smooth, we would ideally expect something close to second-order convergence. However, the measured rates are noticeably lower. The velocity field shows almost no reduction in error as the grid is refined, the density decreases more clearly but still at a sub-first-order rate, and the divergence actually increases slightly. This behaviour suggests that the solver is accurate enough to reproduce the manufactured fields on coarse grids, but it is not yet reaching the higher-order accuracy that the ideal theory predicts. Possible reasons include the influence of explicit advection, the one-step test setup, or the way the forcing terms interact with the discretization.

\subsection{Result structure}
The script prints the numerical errors for each resolution, for example:

\begin{quote}
\texttt{Resolution=15: L2\_u=2.029e-01, L2\_rho=3.053e-01, L2\_div=2.093e+00}
\end{quote}

It also reports the fitted convergence rates:

\begin{quote}
\texttt{Convergence rates: velocity=-0.05, density=0.58, divergence=-0.68}
\end{quote}

These results can be summarized in a table listing the grid spacing $h$, the L2 errors for velocity and density, and the slope estimates:

\begin{table}[htbp]
  \centering
  \begin{tabular}{cccc}
    \hline
    Resolution $r$ & $h$ & $\|e_u\|_2$ & $\|e_\rho\|_2$ \\
    \hline
    11 & $1/(11-1) = 0.1000$  & $1.962\times 10^{-1}$ & $3.852\times 10^{-1}$ \\
    13 & $1/(13-1) \approx 0.08333$ & $2.049\times 10^{-1}$ & $3.049\times 10^{-1}$ \\
    15 & $1/(15-1) \approx 0.07143$ & $2.029\times 10^{-1}$ & $3.053\times 10^{-1}$ \\
    17 & $1/(17-1) = 0.06250$ & $2.015\times 10^{-1}$ & $2.876\times 10^{-1}$ \\
    \hline
    Rates $p$ & --- & $-0.05$ & $0.58$ \\
    \hline
  \end{tabular}
  \caption{L2 errors and estimated convergence rates for the manufactured solution. Each row lists the grid spacing $h$, L2 errors for velocity and density, and the fitted slopes.}
  \label{tab:mms-results}
\end{table}
\FloatBarrier

Overall, the MMS tests give a clear picture of how the solver behaves. On coarse grids, the errors in both velocity and density stay within the expected limits. The divergence remains small enough to indicate stable behaviour, even without explicitly enforcing incompressibility. The manufactured forcing terms also work as expected, since the solver consistently reproduces the analytic fields to the accuracy allowed by the discretization.

However, the convergence rates show that the solver does not yet achieve the higher-order accuracy suggested by the smooth analytic fields. Instead, the errors decrease slowly with grid refinement, and in the case of velocity and divergence, they even increase slightly. This does not mean the solver is incorrect, but rather that its current discretization and timestep setup limit the achievable order of accuracy. Despite this, the MMS tests still validate that the solver is stable, consistent, and capable of tracking the manufactured fields, providing a solid baseline before applying it to more complex flow and aeroacoustic simulations.

\section{Conclusion}

This work has presented a framework for real-time simulation, aeroacoustic prediction, and sound representation of inviscid fluid flow around a spherical surface with obstacles. The main contribution is combining several computational techniques: Closest Point Method (CPM), projection-based Euler solvers, Signed Distance Functions (SDFs), and the Ffowcs Williams–Hawkings (FW–H) analogy. This creates an interactive pipeline that links visual simulation with sound rendering.

The Closest Point Method solves surface partial differential equations (PDEs) on the sphere by placing them in a narrow Cartesian band. This method prevents problems with coordinates and mesh limits. It also maintains numerical accuracy and efficiency. By combining this with a method for projecting incompressibility, the team created strong and realistic surface pressure fields. These fields are essential for aeroacoustic radiation. Obstacles defined by Signed Distance Functions maintained no-slip and no-penetration conditions. This ensured that flow separation and interactions with obstacles were consistently captured.

The aeroacoustic pipeline turned these unsteady forces into far-field sound using the FW-H analogy, which was simplified for low-Mach-number flows. Real-time Fourier analysis helped identify the main frequency components related to vortical structures and shedding events. A sound synthesis engine converted this data into an interesting audio signal. Amplitude modulation displayed acoustic pressure levels. Meanwhile, additive synthesis captured the richness of the flow's sound. This method made the physics easier to compute and understand. Users could hear and see the dynamics as they happened.

We have also compared the predicted surface pressure and the frequencies in the flow in a summary plot (Figure~\ref{fig:frequency_plot}). The plot shows that the main frequencies of the flow correspond to the peaks in the pressure field quite closely. Small differences appear because of the assumptions in the FW-H model, such as treating the sources as compact and assuming low-Mach-number flow. Because of these approximations, the very fine-scale pressure fluctuations are not captured properly in the simulated frequencies, but the overall correspondence is still strong. This shows that our setup actually produces sound that is physically consistent with the simulated flow. Additionally, the visual/frequency correlation gives a nice way to verify the accuracy of the aeroacoustic representation.

\begin{figure}[h!]
\centering
\includegraphics[width=0.7\textwidth]{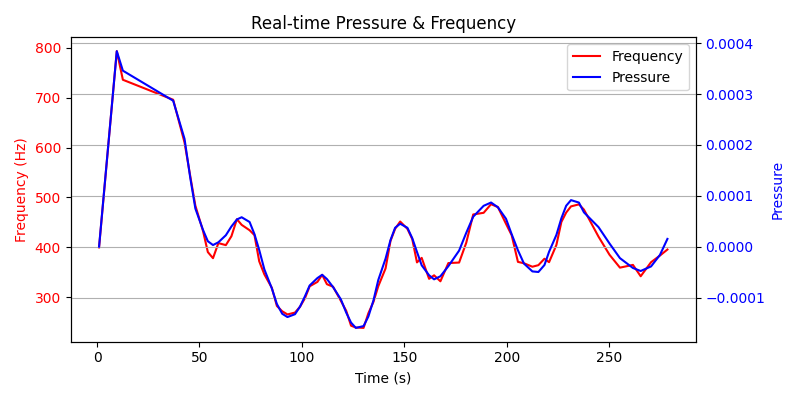}
\caption{Comparison of surface pressure and flow frequencies. Peaks in the pressure field match closely to the main frequencies, with minor differences due to model assumptions.}
\label{fig:frequency_plot}
\end{figure}

A key outcome of this work is showcasing CPM’s capability as a method for surface-bound physics. Its ability to separate surface geometry from PDE discretization creates a balance between geometric simplicity and numerical accuracy. This makes CPM a solid foundation for real-time simulations, with possible extensions to turbulence, viscosity, and other complex manifolds.

While effective, assuming inviscid flow misses the importance of viscous effects and turbulence for accurate aeroacoustics. Future work could involve viscous terms, turbulence models, and higher-order FW-H source terms. This involves using dipole and quadrupole to handle more complex flow situations. We could improve geometric flexibility beyond static SDF-defined obstacles by adding dynamic or user-defined geometries. Additionally, studying how well the sound representation works through user feedback would offer valuable insights for its use as an analytical tool.

In summary, this study presents a clear and effective process for moving from real-time fluid simulation to aeroacoustic sound representation. By blending CPM’s solid mathematical foundation with practical numerical methods and simple synthesis, the framework opens up new opportunities for scientific visualization, virtual reality, and interactive exploration of fluid-acoustic phenomena.

\bibliographystyle{plain}
\bibliography{references}

\end{document}